\documentclass[a4paper]{amsart}

\usepackage[latin1]{inputenc}
\usepackage{amsmath,amssymb,amsthm}

\def\hideJ($J$-){}
\newcommand{\N}{\mathbb{N}}
\newcommand{\R}{\mathbb{R}}
\newcommand{\C}{\mathbb{C}}
\renewcommand{\i}{\mathrm{i}}
\newcommand{\Jets}{\mathcal{J}}
\renewcommand{\S}{\mathfrak{S}}
\newcommand{\typelisse}[1][1]{\Delta^{#1}_{\mathrm{reg}}}
\renewcommand{\binom}[2]{C^{#1}_{#2}}
\newcommand\tsderiv[3][]{{\textstyle{\deriv[#1]{#2}{#3}}}}
\newcommand\deriv[3][]{\frac{\partial^{#1}#2}{\scanDerivVars#3,DerivNoMoreVars,}}
\def\scanDerivVars#1,{%
\def\temp{#1}%
\def\tempLast{DerivNoMoreVars}%
\ifx\temp\tempLast\else%
\partial #1
\expandafter\scanDerivVars\fi%
}
\theoremstyle{plain}
\newtheorem{theorem}{Theorem}
\newtheorem{lemma}{Lemma}
\newtheorem{corollary}[lemma]{Corollary}
\newtheorem{proposition}[lemma]{Proposition}
\theoremstyle{definition}
\newtheorem{definition}[lemma]{Definition}
\theoremstyle{remark}
\newtheorem{remark}{Remark}

\begin{document}

\title[Type of real hyper-surfaces]{Regular type of real hyper-surfaces
  in (almost) complex manifolds}

\author[J.-F.~Barraud]{Jean-François Barraud}%
\address{UMR 8524, Université Lille 1\\
59 655 Villeneuve d'Ascq Cedex - FRANCE}
\author[E.~Mazzilli]{Emmanuel~Mazzilli}%
\address{UMR 8524, Université Lille 1\\
59 655 Villeneuve d'Ascq Cedex - FRANCE}


\begin{abstract}
  The regular type of a real hyper-surface $M$ in an (almost) complex manifold
  at some point $p$ is the maximal contact order at $p$ of $M$ with
  germs of non singular (pseudo) holomorphic disks. The main purpose of
  this paper is to give two intrinsic characterizations the type~: 
  one in terms of Lie brackets of a complex tangent vector field on $M$,
  the other in terms of some kind of derivatives of the Levi form.
\end{abstract}

\maketitle

The purpose of this paper is to study order of contact between
holomorphic curves (as well pseudoholomorphic curves) and a real
hyper-surface. In the integrable case, this invariant is connected to the
boundary behavior of the Cauchy-Riemann equations, the Bergman kernel,
invariant metrics, etc.. see for example \cite{Da1} \cite{Ca1, Ca2},
\cite{Mc}, \cite{Kho}. One of the difficulty to use this number in
function theory on a domain $D$, is due to the fact that in general, it
was not known how to calculate it with intrinsic ``complex geometry'' of
the boundary of $D$. In $\Bbb{C}^2$ the situation is clear (see
\cite{Kho}) ; for $\Bbb{C}^n$ ($n>2$), we only know how to compute the
order of contact of complex hyper-surface with the natural Lie algebra of
the boundary of $D$ (see \cite{Blo-Gra}). In \cite{Blo-Gra}, I.~Graham
and T.~Bloom ask how to characterize the regular ``type'' in a similar
intrinsic way with only one complex vector field; T.~Bloom
(see\cite{Blo}) succeeds in the pseudoconvex case in $\Bbb{C}^3$ but
unfortunately, the result is not valid without the pseudoconvexity
hypothesis (we think that even with the pseudoconvexity, this result is
false in $\Bbb{C}^n$ with $n\geq4$).

In this article, we consider an hyper-surface (pseudoconvex or not) in
$\C^{n}$ or in $\R^{2n}$ endowed with an almost complex structure, and
characterize its ``regular one type'', by means of Lie brackets of one
``complex tangent vector field'' on the hyper-surface (see theorem $1$),
and by means of derivatives of the ``levi form'' in some way (see theorem
2). We decided to work in the non-integrable situation, because our main
argument shows up naturally in this setting. In particular, we recall the
definition of the Levi form in the non-integrable case, show that it has
the same properties as the classical one, and use it to define strictly
pseudoconvexity and pseudoconvexity.

The plane of the paper is the following: in the first section, we define
relevant objects in the non-integrable case and we expound the main
results; in section two, we compare jets of holomorphic curves (with
contact order $k$) with jets of complex tangent vector fields;
in section three, we relate the order of contact with complex line
sub-bundles of $T^{J}(M)=TM\cap JTM$ (see the definition in the first
section) nearly involutive (see theorem 1); in section four, we introduce
the ``derivatives of the Levi form'' and we demonstrate the theorem 2; in
the last section, we discuss the base point dependency of the type in
(complex) dimension 2 and in a particular case in dimension 3.

\section{Notations and main results}

Let us first recall, in the almost complex case, the notions of ``complex
tangent space'', ``complex vector fields'', ``Levi form'' and
``pseudoconvexity'' and derive some of their properties well known in the
integrable case.

\medskip

All along this paper, we will work with the following objects~: let $J$
be an almost complex structure on $\R^{2n}$ (with $J(0)=\i$) and $M$ an
oriented real hyper-surface through $0$ in $\R^{2n}$. Let $T^{J}M=TM\cap
JTM$ denote the $J$-invariant part of $TM$. The sections of $T^{J}M$ will
be called $J$-tangent or complex tangent vector fields on $M$. We choose
a vector field $N$ transversal to $M$ so that we have~: 
$$
\R^{2n}=TM\oplus\R N \quad\text{and}\quad TM=T^{J}M\oplus\R(JN).
$$
Let $\phi~:\R^{2n}\rightarrow\R$ be a defining function for $M$, {\it
  i.e.} a function such that $M=\phi^{-1}(0)$, $0$ is a regular value of
$\phi$, and $\phi$ defines the positive orientation of $M$ .  Of course,
if $\psi$ is another defining function of $M$, then there exists a non
negative function $f:\R^{2n}\rightarrow]0,+\infty[$ such that
$\psi=f\phi$. Once a defining function $\phi$ and a metric on $\R^{2n}$
are chosen, we can take $N$ to be the gradient $\nabla\phi$ of $\phi$
with respect to the metric.

The set $M_{+}=\{\phi>0\}$ will be called the outside, and
$M_{-}=\{\phi<0\}$ the inside defined by $M$.

The main purpose of this paper is to characterize the maximal contact
order $\typelisse[1](p)$ of $M$ with regular (eventually pseudo-)
holomorphic disks at some point $p$ intrinsically, {\it i.e.} in terms of
tangent vector fields on $M$. Recall that a germ of map
$u~:(\C,0)\rightarrow(\R^{2n},0)$ is pseudo-holomorphic with respect to
$J$, or $J$-holomorphic if it satisfies the relation $du+J(u)du\i=0$.
Since this should not produce any ambiguity, we will speak about
``holomorphic'' disks in both the integrable and the non integrable case,
forgetting about the ``pseudo'' or $J$ prefix. The precise notions of
contact and type we will work with are the following~:
\begin{definition}
  Let $u:(\Delta,0)\rightarrow(\R^{2n},0)$ be a (pseudo-)holomorphic disk
  mapping $0$ to $0$, and regular {i.e.} such that $du(0)\neq0$. Its
  contact order $\delta_{0}(M,u)$ with $M$ at $0$ is the degree of the
  first term in the Taylor expansion of $\phi\circ u$.
\end{definition}

\begin{definition}
  Let $\mathcal{D}_{\mathrm reg}$ be the set of all germs of
  (pseudo-)holomorphic disks $u$ such that $u(0)=0$ and $du(0)\neq 0$. The
  (regular) type of $M$ at $0$ is the number
  $$
  \typelisse[1](M,0)=\sup_{u\in\mathcal{D}_{\mathrm{reg}}} (\delta_{0}(M,u))
   $$
\end{definition}

Our aim is to compute  the type of $M$ at $0$ in terms of tangent vector
fields on $M$. We propose two results in this direction, inspired from
two classical results in complex analysis. 

Our first  result in this direction is a generalization of the classical
results of \cite{Kho} and \cite{Blo-Gra} and computes the type in term of
Lie brackets of ``one'' tangent vector field.
\begin{theorem}\label{thm:Contact order via Lie brackets}
  Let $\delta=\C X$ be $J$-line in $T_{0}M$. The three following
  conditions are equivalent~:
  \begin{enumerate}
  \item\label{item:Crochets de Lie-u} 
    There exists a smooth ($J$-)holomorphic disk $C=u(\Delta)$ tangent to
    $M$ at $0$ with order $k+2$, with direction $\delta$ ({\it i.e.
      $T_{0}C=\delta$}).
  \item\label{item:Crochets de Lie-Droite stable} 
    There exists a complex line sub-bundle $L$ of $TM$ such that
    $L_{0}=\delta$ and all the Lie brackets of length at most $k+1$ of
    sections of $L$ steel belong to $\delta$ at $0$.
  \item\label{item:Crochets de Lie-X commute} 
    there exists a $J$-tangent vector field $X_{0}$ with
    $X_{0}(0)\in\delta$ such that all the Lie brackets of $X_{0}$ and
    $JX_{0}$ of length at most $k+1$ vanish at $0$.
  \end{enumerate}
  Moreover, if $\nabla$ is a symmetric connexion, $u$ and $X_{0}$ can be
  chosen such that at $0$~: 
  $X_{0}=\deriv{u}{x},\nabla_{X_{0}}X_{0}=\deriv[2]{u}{x^{2}}, \dots,
  \nabla_{X_{0}}(\nabla_{X_{0}}(\dots \nabla_{X_{0}}X_{0}))=
  \deriv[k+1]{u}{x^{k+1}}$.
\end{theorem}

\begin{remark}
 If $L$ is involutive, i.e. (every where~!) stable by Lie brackets, then by
 the Frobenius theorem, $M$ is foliated by $J$-holomorphic curves.
\end{remark}

Our second result in the same direction uses ``higher order'' Levi forms,
which are some kind of derivatives of the standard Levi form~; in
contrast to the standard Levi form, they depend on the choice of a
symmetric connexion.

Let us first discuss the Levi form itself. It is often defined on the
complexification of $TM$, but we choose not to work in this setting to
reduce the number of (almost) complex structures involved.

Thus we will use the following definition~:
\begin{definition}
  Let $T^{J}M=TM\cap J(TM)$ be the $J$-invariant part of $TM$. The
  function $L_{\phi}:T^{J}M\rightarrow\R$ defined by
  \begin{equation}
  \label{eq:Levi form}
  L_{\phi}(X)=d\phi(J[X,JX])
  \end{equation}
  is in fact a quadratic form on $TM$, called the Levi form associated to
  $\phi$. It depends on $\phi$ only up to multiplication by a non 
  negative function.
\end{definition}

The fact that $L(X)$ does not depend on the derivatives of $X$ is a
straight forward computation. Let us just mention by the way the polar
form associated to $L_{\phi}$~:
  \begin{equation}
  \label{eq:Hermitian Levi form}
  \Theta_{\phi}(X,Y)=\frac{1}{2}\ \big(d\phi(J[X,JY]+J[Y,JX])
                     +\i\ d\phi(J[X,Y]+J[JX,JY])\big)
  \end{equation}
In particular, if $\alpha$ and $\beta$ are two real functions on $M$, then
\begin{equation}
  \label{eq:L((a+bJ)X)}
  L_{\phi}((\alpha+\beta J)X)=(\alpha^{2}+\beta^{2})L_{\phi}(X).
\end{equation}

The classical relation between the Levi form and the second derivative of
$\phi$ has to be modified in the non integrable case to take derivatives
of the almost complex structure $J$ into account~:
\begin{proposition}
  Let $D^{2}_{\nabla}\phi$ denote the second derivative of $\phi$ with
  respect to $\nabla$. Then
  \begin{equation}
    \label{eq:L(X)=h(X,X)+h(JX,JX)+new term}
    L_{\phi}(X)=D^{2}_{\nabla}\phi(X,X)+D^{2}_{\nabla}\phi(JX,JX)
    +d\phi\Big((\nabla_{JX}J)\,X-\nabla_{X}J)\,JX\Big)
  \end{equation}
\end{proposition}

\begin{remark}
  We will use a multiplicative notation for the derivatives, i.e. when
  differentiating many times in different directions,
  $\nabla_{X_{1}}(\nabla_{X_{2}}(\dots\nabla_{X_{k-1}}X_{k}))$ will be
  replaced by $X_{1}\cdot X_{2}\cdots X_{k}$, and if all the $X_{i}$ are
  the same, simply by $X^{k}$. We stress that this notation is somewhat
  misleading since this ``product'' is not associative.
\end{remark}

The particular cases when the Levi form is non negative or positive are
specially interesting.
\begin{definition}
  The hyper-surface $M$ is said to be
  \begin{itemize}
  \item          pseudoconvex if $L_{\phi}\geq0$.
  \item strictly pseudoconvex if $L_{\phi}>0$.
  \end{itemize}
\end{definition}

The Levi form is a good tool to study the holomorphic disks tangent to
$M$ at the first order, i.e. to check whether the type is $2$ or grater than
$2$~:

\begin{proposition}\label{thm:L>0 et courbe d'appui}
  Consider the Levi form at the origin.
  \begin{itemize}
  \item
    If $\exists X\in T^{J}_{0}M,\ L_{\phi}(X)>0$, then there exits a
    smooth \hideJ($J$-)holomorphic disk $u:\Delta\rightarrow\R^{2n}$ tangent to
    $M$ at $0$ from the outside (i.e. lying in $M_{+}$). \item
    $\exists X\in T^{J}_{0}M,\ L_{\phi}(X)=0$ if and only if there exits
    a smooth \hideJ($J$-)holomorphic disk $u:\Delta\rightarrow\R^{2n}$
    tangent to $M$ at $0$ with contact order $>2$.
  \end{itemize}
 Moreover, in each case, $X$ and $u$ can be chosen such that at $0$~:
 $X=\deriv{u}{x}$.
\end{proposition}
\begin{remark}
This is just a computation, which will be detailed later on, as an
introduction to the definition of higher order Levi forms. The existence of such holomorphic disk 
is useful to study properties of hyperbolicity for strictly pseudoconvex hyper-surfaces (see \cite{Deb}).  
\end{remark}
The main result is as follows~:

\begin{theorem}\label{thm:Contact order via higher Levi forms}
  Out of the derivatives of the Levi form (see definition \ref{def:Lpq} in section \ref{sec:Higher Levi Forms} for precise statement),
  one can compute functions 
  $L^{p,q}: (\R^{2n})^{p+q+1}\rightarrow\R$ for all $p,q\in\N$ such that
  the two following conditions are equivalent~:
  \begin{enumerate}
  \item \label{item:Contact via Levi form-u exists}
    There exists a regular ($J$-)holomorphic disk tangent to $M$ at $0$
    with order $k+2$ 
  \item \label{item:Contact via Levi form-Lpq(X)=0}
    There exists a complex tangent vector field $X\in\Gamma(T^{J}M)$ such
    that
    $$
    \forall (p,q), p+q\leq k-1,\ 
    L^{p,q}(X(0),X^{2}(0),\dots,X^{p+q+1}(0))=0.
    $$
  \end{enumerate}
  Moreover, $u$ and $X$ can be chosen such taht at $0$~:
  $X=\deriv{u}{x},\dots ,X^{k+1}=\deriv[k]{u}{x^{k+1}}$~; 
(recall that $X^{p}$ stands for
$\underset{p}{\underbrace{\nabla_{X}(\cdots(\nabla_{X}(X)))}}$).
\end{theorem}

\begin{remark}
This result generalizes the theorem of Bloom in $\Bbb{C}^3$ (see \cite
{Blo}) for pseudoconvex hyper-surface; but compare carefuly the definition
of iterated derivatives of the levi form in our theorem and in \cite
{Blo}.
\end{remark}

The two theorems (\ref{thm:Contact order via Lie brackets} and
\ref{thm:Contact order via higher Levi forms})  rest on the comparison of
jets of \hideJ($J$-)holomorphic disks and the appropriate notion of jets
of vector fields on $M$.

\section{Jets of disks and tangent vector fields}

\begin{definition}
  Define the $(p,q)$-th derivative of some vector field $X$ in $\R^{2n}$,
  to be the vector $X$ differentiated $p$ times in its own direction and
  then $q$ times in the $JX$ direction~:
  $$ 
  D^{p,q}_{X}X(0)=(JX)^{q}X^{p}\cdot X(0)=
  \underset{q}{\underbrace{(JX\cdot \cdots (JX}}\cdot 
  \underset{p}{\underbrace{( X\cdot \cdots ( X}}\cdot X)))
  (0)
  $$
\end{definition}

\begin{definition}
  The collection of all the $(p,q)$-th derivatives of a vector field $X$
  on $\R^{2n}$ at $0$ for $0\leq p+q\leq k$ will be called its oder $k$
  jet and denoted by $j^{k}_{0}X$. This defines a map 
  $$
  j^{k}_{0}:
  \begin{array}{rcl}
    \Gamma(T\R^{2n})&\rightarrow &\Jets^{k}_{0}(\R^{2n})=(\R^{2n})^{\frac{(k+1)(k+2)}{2}}\\
    X&\mapsto&(X(0),\dots, D^{p,q}_{X}X(0),\dots, D^{0,k}_{X}X(0))
  \end{array}
  $$
  and a jet $\xi\in\Jets^{k}_{0}(\R^{2n})$ is said to be
  \emph{realizable} if it belongs to
  $\Jets^{k}_{0}(M)=j^{k}_{0}(\Gamma(T^{J}M))$.
\end{definition}

Given an holomorphic disk $u$ tangent to $M$, one would like to know
whether or not it's jet at $0$ can be realized as the jet of some vector
field $X\in\Gamma(T^{J}M)$. This is what this section is devoted to.

Let us start with a basic observation~:
\begin{proposition}
  Let $X$ and $Y$ be two $J$-tangent vector fields such that
  $j^{k}_{0}(X)=j^{k}_{0}(Y)$. Then for all $p,q\in\N$ with $p+q=k+1$, we
  have
  $$
  D^{p,q}_{X}X(0)-D^{p,q}_{Y}Y(0)\in T^{J}M
  $$
\end{proposition}
\begin{proof}
  One has to prove the two equalities
  \begin{align}
    d\phi(D^{p,q}_{X}X(0)-D^{p,q}_{Y}Y(0))&=0
    \label{eq:DpqX-DpqY in TM}\\
    d\phi(J[D^{p,q}_{X}X(0)-D^{p,q}_{Y}Y(0)])&=0.
    \label{eq:DpqX-DpqY in TJM}
  \end{align}
  Let us compute $D^{p,q}_{X}(d\phi(X))$. On one hand, it is $0$ (since
  $X\in TM$). On the other hand, the development of
  $D^{p,q}_{X}(d\phi(X))$ is a sum of terms of the form
  $$
  D^{\alpha}\phi(D^{p_{1},q_{1}}_{X}X(0),\dots,
  D^{p_{\alpha},q_{\alpha}}_{X}X(0))
  $$
  If $p_{i}+q_{i}\leq k$, then this term is the same as
  $$
  D^{\alpha}\phi(D^{p_{1},q_{1}}_{Y}Y(0),\dots,
  D^{p_{\alpha},q_{\alpha}}_{Y}Y(0))
  $$
Thus $0-0=D^{p,q}_{X}(d\phi(X))-D^{p,q}_{Y}(d\phi(Y)) =
d\phi(D^{p,q}_{X}X-D^{p,q}_{Y}Y)$.

The equality (\ref{eq:DpqX-DpqY in TJM}) is obtained along the same lines
(eventually taking the derivatives of $J$ into account) using
$D^{p,q}_{X}(d\phi(JX))=0$ since $JX\in TM$.
\end{proof}

Realizability of jets can then be recursively tested as follows~:

\begin{proposition}
  \label{thm:Jets realisables}
  Let $\xi\in \Jets^{k+1}_{0}(\R^{2n})$. Let $[\xi]_{k}$ be its
  $\Jets^{k}_{0}(\R^{2n})$-component, and
  $(\xi_{k+1,0},\dots,\xi_{p,q},\dots,\xi_{0,k+1})$ its homogeneous part
  of degree $k+1$. Then
  \begin{equation}
   \exists X\in\Gamma(T^{J}M), j^{k+1}_{0}(X)=\xi\label{eq:xi realisable}
  \end{equation}
if and only if
  \begin{equation}\label{eq:xi realisable ordre precedent et compatible}
    \left\{
    \begin{aligned}
    &\exists X_{1}\in\Gamma(T^{J}M), j^{k}_{0}(X_{1})=[\xi]_{k}\\
    &\forall (p,q), p+q=k+1,\ \ 
      \xi_{p,q}-[j^{k+1}_{0}(X_{1})]_{p,q}\in T^{J}_{0}M
    \end{aligned}\right.
  \end{equation}
\end{proposition}

In other words, if you take a realizable $k+1$-jet, then the $N$ and $JN$
components of its maximal degree part are completely determined by its
$k$-order part, but the $T^{J}M$ component can be anything.

\begin{proof}
(\ref{eq:xi realisable})$\Rightarrow$(\ref{eq:xi realisable ordre
  precedent et compatible}) is obvious. So let us start with $X_{1}$ and
the compatibility condition (\ref{eq:xi realisable ordre precedent et
  compatible}). Let $(T_{1},\dots,T_{n-1})$ be a (complex) basis of
$T^{J}M$ and let $X=X_{1}+\sum \lambda_{i}T_{i}$, where the $\lambda_{i}$
are complex valued functions on $M$, all the derivatives of which vanish
at $0$ up to order $k$, and such that $\lambda_{i}(0)=1$. Then if
$p+q\leq k$, $D^{p,q}_{X}X=D_{X_{1}+\sum\lambda_{i}T_{i}}(X+\sum\lambda
_{i}T_{i})$ is the sum of $D^{p,q}_{X_{1}}X_{1}$ and of terms involving
derivatives of the $\lambda_{i}$ of order at most $k$. Therefore at $0$
we have 
\begin{equation}
\label{eq:DpqX(0)=DpqX1(0)}
\forall(p,q),p+q\leq k,\quad D^{p,q}_{X}X(0)=D^{p,q}_{X_{1}}X_{1}(0)
\end{equation}

For $p+q=k+1$, the $k+1$-th derivatives of $\lambda_{i}$ appear in the
following terms~:
\begin{align}
  D^{p,q}_{X}X(0)
  &=D^{p,q}_{X_{1}}X_{1}(0)+\sum(D^{p,q}_{X_{1}}\lambda_{i})(0)T_{i}(0)
  \notag\\
  &=D^{p,q}_{X_{1}}X_{1}(0)+
  \sum D^{k+1}_{0}\lambda_{i}(JX_{1},\dots,JX_{1},X_{1},\dots,X_{1})T_{i}(0)
\end{align}
which proves that the $k+1$ degree part of $j^{k}_{0}(X_{1})$ can only be
modified in the $T^{J}M$ direction, and all the modification in this
direction are possible since any
$(\alpha_{k+1,0},\alpha_{k,1},\dots,\alpha_{0,k+1})\in\C^{k+2}$ can be
realized as the derivatives of a complex valued function in the
directions $X_{1}$ and $JX_{1}$ (provided that $X_{1}\neq 0$)~:
$\forall(\alpha_{k+1,0},\alpha_{k,1},\dots,\alpha_{0,k+1})\in\C^{k+2}
\quad\exists\lambda: M\rightarrow\C$ such that
$$
\left\{
  \begin{aligned}
    &\lambda(0)=1\\
    &\forall r\leq k,\quad D^{r}_{0}\lambda=0\\
    &\forall(p,q),p+q=k+1, \quad
       D^{k+1}_{0}\lambda(JX_{1}\dots JX_{1},X_{1}\dots X_{1})=\alpha_{p,q}
  \end{aligned}
\right.
$$
\end{proof}

\begin{corollary}\label{thm:u<->X}
  There exists a regular holomorphic disk $u:\Delta\rightarrow\R^{2n}$
  tangent to $M$ at $0$ with order $k+2$ ({\it i.e.} $\phi\circ
  u(z)=O(z^{k+2})$) if and only if the $k+1$-jet of $u$ at $0$ is
  realizable, {\it i.e.}
  $$
  \exists X\in\Gamma(T^{J}M),\quad
  \forall(p,q),p+q\leq k,\quad
  D^{p,q}_{X}X(0)
  =\deriv[p+q]{u}{x^{p},y^{q}}\,\deriv{u}{x}(0).
  $$
\end{corollary}

\begin{proof}
  Let $u$ be a disk tangent to $M$ with order $k+3$, $\xi$ the $k$-jet of
  $\deriv{u}{x}$ at $0$ ({\it i.e.} $\xi_{p,q}=
  \deriv[p+q+1]{u}{x^{p+1},y^{q}}(0)$). By induction, there exists a $J$
  tangent vector field $X_{1}$ realizing the $k$-jet of $\deriv{u}{x}$.
  Then, comparing the developments of 
  $$
  \deriv[k+1]{}{x^{p},y^{q}}\Big[d\phi({\tsderiv{u}{x}})\Big](0)
  \qquad\text{ and }\qquad
  D^{p,q}_{X_{1}}\Big[d\phi(X_{1})\Big](0)
  $$
  for $p+q=k+1$, it turns out that all the terms are the same except
  maybe the terms where no derivative is applied to $d\phi$. Thus, we
  have~:
  \begin{equation}
    \label{eq:pq-deriv dphi(u')=pq-deriv dphi(X)}
  \tsderiv[k+1]{}{x^{p},y^{q}}\Big[d\phi(\tsderiv{u}{x})(0)\Big](0)-
  D^{p,q}_{X_{1}}\Big[d\phi(X_{1})\Big](0)=
  d\phi\Big[\tsderiv[k+1]{}{x^{p},y^{q}}\,{\tsderiv{u}{x}(0)}-
  D^{p,q}_{X_{1}}X_{1}(0)\Big].
  \end{equation}
  On the other hand, $\deriv[k+1]{}{x^{p},y^{q}}
  \Big[d\phi({\tsderiv{u}{x}})\Big](0)=0$ since $\phi\circ
  u(z)=O(z^{k+3})$, and $D^{p,q}_{X_{1}}\Big[d\phi(X_{1})\Big](0)=0$
  since $X_{1}\in TM$. Thus we have
  \begin{equation}
    \label{eq:N-composante de u' et X1}
    d\phi\Big(\xi_{p,q}-[j^{k+1}_{0}X_{1}]_{p,q}\Big)=0.
  \end{equation}
  
  Comparing in the same way
  $$
  \deriv[k+1]{}{x^{p},y^{q}}\Big[d\phi(J(u){\tsderiv{u}{x}})\Big](0)
  \qquad\text{ and }\qquad
  D^{p,q}_{X_{1}}\Big[d\phi(JX_{1})\Big](0)
  $$
  we get 
  \begin{equation}
    \label{eq:pq-deriv dphi(Ju')=pq-deriv dphi(JX)}
  \tsderiv[k+1]{}{x^{p},y^{q}}\Big[d\phi(\tsderiv{u}{y})(0)\Big](0)-
  D^{p,q}_{X_{1}}\Big[d\phi(JX_{1})\Big](0)=
  d\phi\Big[J\tsderiv[k+1]{}{x^{p},y^{q}}\,{\tsderiv{u}{x}(0)}-
  JD^{p,q}_{X_{1}}X_{1}(0)\Big].
  \end{equation}
  and finally
  \begin{equation}
    \label{eq:JN-composante de u' et X1}
    d\phi\Big(J\xi_{p,q}-J[j^{k+1}_{0}X_{1}]_{p,q}\Big)=0.
  \end{equation}
  The relations (\ref{eq:N-composante de u' et X1}) and
  (\ref{eq:JN-composante de u' et X1}) and the proposition \ref{thm:Jets
    realisables} now proves that $X_{1}$ can be modified into a new
  complex tangent vector field $X$ to have~:
  $$
  j^{k+1}_{0}X=\xi.
  $$

  \medskip
  Start now with $u$ and $X\in T^{J}M$ such that $\forall(p,q),p+q\leq k,
  \deriv[p+q+1]{u}{x^{p+1}y^{q}}(0)=D^{p,q}_{X}X(0)$ and compute
  $\deriv[p+q]{\phi\circ u}{x^{p},y^{q}}(0)$ for $p+q\leq k+1$. The main
  tool is again the relations (\ref{eq:pq-deriv dphi(u')=pq-deriv
    dphi(X)}) and (\ref{eq:pq-deriv dphi(Ju')=pq-deriv dphi(JX)}). If
  $p>1$ we have
  \begin{align}
    \deriv[p+q]{}{x^{p},y^{q}}\big[\phi\circ u\big](0)
    &=\deriv[p+q]{}{x^{p-1},y^{q}}
      \big[d\phi({\tsderiv{u}{x}})\big](0)\notag\\
    &=D^{p-1,q}_{X}\big[d\phi(X)\big](0)=0
  \end{align}
  and if $p=0$~:
  \begin{align}
    \deriv[q]{}{y^{q}}\big[\phi\circ u\big](0)
    &=\deriv[q-1]{}{y^{q-1}}
      \big[d\phi(J{\tsderiv{u}{x}})\big](0)\notag\\
    &=D^{0,q-1}_{X}\big[d\phi(JX)\big](0)=0
  \end{align}
Finally, we have $\phi\circ u(z)=O(z^{k+2})$.
\end{proof}

\section{Regular type and Lie algebra of line bundles}

\subsection{Commutativity of vector fields at one point}
\begin{definition}
  A complex tangent vector field $X$ is said to commute at $0$ up to
  order $k$ if and only if $\forall m\leq k$,$\forall
  X_{1},\dots,X_{m}\in\{X,JX\}$ and $\forall\sigma\in\S_{m}$ we have
  $X_{\sigma(1)}\cdot X_{\sigma(2)}\cdots X_{\sigma(k)}=X_{1}\cdot
  X_{2}\cdots X_{k}$.
\end{definition}

Of course, we can use the group structure of $\S_{k}$ and the properties
of the Lie brackets to reduce attention only to a few permutations~:
\begin{proposition}
  Let $X$ be a complex tangent vector fields. The four following
  conditions are equivalent~:
  \begin{enumerate}
  \item\label{item:X commute}
    $X$ commutes up to order $k$ at $0$.
  \item\label{item:X pre-commute}
    $\forall p,q/ p+q\leq k,$
    $$
    \underset{q  }{\underbrace{JX\cdots JX}}
    \cdot\underset{p  }{\underbrace{ X\cdots  X}}
    \cdot JX =
    \underset{q+1}{\underbrace{JX\cdots JX}}
    \cdot\underset{p  }{\underbrace{ X\cdots  X}}
    $$
  \item\label{item:Crochets de Lie=0} 
    All the Lie brackets of $X$ and $JX$ up to length $k$ vanish at $0$.
  \item\label{item:deriv[X,JX]=0}
    $\forall m\leq k-2,\forall X_{1}\dots X_{m}\in\{X,JX\}$,\ \
    $X_{1}\cdots X_{m}\cdot[X,JX]=0$
\end{enumerate}
\end{proposition}

\begin{proof}
  $(\ref{item:X commute})\Rightarrow(\ref{item:X pre-commute})$ is clear.

  $(\ref{item:X pre-commute})\Rightarrow(\ref{item:X commute})$ rests on
  the remark that if $X$ commutes up to order $k-1$, then the relations
  of $(\ref{item:X commute})$ for $m=k$ and $\sigma(k)=k$ are automatic~: 
  $$
  X_{1}\cdots X_{k-1}\cdot X_{k}=
  \sum_{\substack{
      1\leq\alpha\leq k-1\\
      I_{1}\cup\dots\cup I_{\alpha}=\{1\dots n\}\\
      I_{i}\cap I_{j}=\emptyset}}
   D^{\alpha}X_{k}(X^{I_{1}},\dots,X^{I_{\alpha}})
  $$ 
  where $X^{I_{i}}=X_{i_{1}}\cdots X_{i_{r}}$ when
  $I_{i}=\{i_{1}<\dots<i_{r}\}$.

  The group $\S_{k-1}$ acts transitively on the indexing set, and since
  $X$ commutes up to order $k-1$, $X^{I_{i}}$ does not depend on the
  order used on $I_{i}$. We conclude that $X_{1}\cdots
  X_{k}=X_{\sigma(1)}\cdots X_{\sigma(k-1)}\cdot X_{k}$ for all
  $\sigma\in\S_{k-1}$. This ends the proof since the permutations
  corresponding to the relation (\ref{item:X pre-commute}) and $\S_{k-1}$
  generate $\S_{k}$. 

  $(\ref{item:X commute})\Rightarrow(\ref{item:Crochets de Lie=0})$ is
  clear.
  
  $(\ref{item:Crochets de Lie=0}) \Rightarrow
  (\ref{item:deriv[X,JX]=0})$~: suppose that all the Lie brackets of $X$
  and $JX$ vanish at $0$ up to length $k$, and, by induction, that all
  the derivatives of $[X,JX]$ in the $X$ and $JX$ directions vanish at
  $0$ up to order $k-3$. Then in the development of $[X_{1}[\dots
  X_{k-2}[X,JX]]]$, all the terms vanish but the term where $[X,JX]$ is
  differentiated $k-2$ times~:
  $$
   0=[X_{1}[\dots X_{k-2}[X,JX]]]=X_{1}\cdots X_{k-2}\cdot[X,JX]
  $$
  
  $(\ref{item:deriv[X,JX]=0}) \Rightarrow (\ref{item:X commute})$~: it is
  enough to prove the relation when $\sigma$ is a transposition
  $(i,i+1)$. But $X_{1}\cdots X\cdot (JX)\cdots X_{m}-X_{1}\cdots
  (JX)\cdot X\cdots X_{m} = X_{1}\cdots[X,JX]\cdots X_{m}$, and all the
  terms in the development of this last expression involve a derivative
  of order at most $m-2$ of $[X,JX]$ in the $X$ and $JX$ directions~:
  thus it is null.
\end{proof}

We can now come to the proof of theorem \ref{thm:Contact order via Lie
  brackets}.

\subsection{Proof of theorem \ref{thm:Contact order via Lie
  brackets}}

\paragraph{Proof of (\ref{item:Crochets de Lie-u}) $\Rightarrow$
  (\ref{item:Crochets de Lie-X commute}) in theorem \ref{thm:Contact
    order via Lie
  brackets}~:}

Start with a regular disk $u:(\Delta,0)\rightarrow(\R^{2n},0)$ tangent to
$M$ with order $k+2$. Then, by corollary \ref{thm:u<->X}, there exists a
complex tangent vector field $X$ such that $\forall (p,q),p+q\leq k,
\deriv[p+q]{}{x^{p},y^{q}}\deriv{u}{x}(0)=D^{p,q}_{X}X(0)$. Let us check
that this $X$ commutes up to order $k+1$ at $0$, {\it i.e.} that
$D^{p,q}_{X}(JX)(0)=D^{p-1,q+1}_{X}X(0)$ for $p+q\leq k$.

We have the relation
$$
 D^{p,q}_{X}(JX)=\sum_{\substack{a+c=p\\b+d=q}}
 \binom{a}{p}\binom{b}{q}\ (D^{a,b}_{X}J)\ D^{c,d}_{X}X
$$
Differentiating $\deriv{u}{y}=J(u)\deriv{u}{x}$ with respect to $x$ $p$
times and to $y$ $q$ times, we obtain
$$
 \tsderiv[p+q]{}{x^{p},y^{q}}\,\tsderiv{u}{y} 
=\sum_{\substack{a+c=p\\b+d=q}} \binom{a}{p}\binom{b}{q}
 \ (\tsderiv[a+b]{}{x^{a},y^{b}}(J(u)))
 \ \tsderiv[c+d]{}{x^{c},y^{d}}\tsderiv{u}{x}
$$
fom which we deduce that~:
\begin{align*}
D^{p,q}_{X}(JX)(0)
&=\deriv[p+q]{}{x^{p},y^{q}}\,\tsderiv{u}{y}(0) \\
&=\deriv[p+q]{}{x^{p-1},y^{q+1}}\,\tsderiv{u}{x}(0) \\
&=D^{p-1,q+1}_{X}X(0)
\end{align*}

This ends the proof of (\ref{item:Crochets de
  Lie-u})$\Rightarrow$(\ref{item:Crochets de Lie-X commute}) in theorem
\ref{thm:Contact order via Lie brackets}.

\medskip

\paragraph{Proof of (\ref{item:Crochets de Lie-X commute})
  $\Rightarrow$(\ref{item:Crochets de Lie-u}) in theorem \ref{thm:Contact
  order via Lie brackets}~:}

Let now $X$ be a vector field commuting up to order $k+1$ at $0$. By an
argument of J.-C.~Sikorav \cite{Sik}, there exists a germ
$u:\Delta\rightarrow\R^{2n}$ of $J$-holomorphic disk such that $\forall
m\leq k+1$, $\deriv[m]{u}{x^{m}}(0)=X^{m}(0)$. Since $X$ commutes, one
can easily check that
$\deriv[p+q]{}{x^{p},y^{q}}\,\deriv{u}{x}=D^{p,q}_{X}X(0)$ and then apply
proposition \ref{thm:u<->X}.

\paragraph{Proof of (\ref{item:Crochets de Lie-Droite stable})
  $\Leftrightarrow$(\ref{item:Crochets de Lie-X commute}) in theorem
  \ref{thm:Contact order via Lie brackets}~:}

The basic remark here is that if $X$ is replaced by $(\alpha+\beta J)X$
where $\alpha$ and $\beta$ are real valued functions on $M$, then an
iterated Lie bracket of $X$ and $JX$ of length $k$ is affected by
addition of a linear combination of strictly shorter Lie brackets.
Therefore, if $X$ commutes up to roder $k+1$ at $0$, then the line
subbundle generated by $X$ is stable at $0$ up to order $k+1$. In the
other direction, if a line bundle is stable at $0$ up to order $k+1$,
then choosing a local non vanishing section $X$ of it, and suitable
functions $\alpha$ and $\beta$, we obtain a vector field commuting at $0$
up to order $k+1$.



\section{Higher order Levi forms}
\subsection{Definition~:}\label{sec:Higher Levi Forms}
As an introduction, let us discuss the proof of proposition \ref{thm:L>0
  et courbe d'appui}. The main tool here is the relationship between the
Levi form and laplacian of $\phi\circ u$ for tangent disks $u$.

Suppose that $X\in T^{J}_{0}M$ satisfies $L_{\phi}(X)>0$. Recall that
from an argument of Sikorav \cite{Sik}, all the derivatives
$\deriv[k]{u}{x^{k}}(0)$ of a \hideJ($J$-)holomorphic disk can be chosen
arbitrarily, {\it i.e.} for all $X_{1},\dots,X_{k}\in\R^{2n}$ there
exists a germ of holomorphic disk such that
$X_{i}=\deriv[i]{u}{x^{i}}(0)$.

Let $u:\Delta\rightarrow\R^{2n}$ be a (germ of) \hideJ($J$-)holomorphic disk such
that $u(0)=0$, $\deriv{u}{x}(0)=X$. In the Taylor expansion of $\phi\circ
u$, use the fact that $u$ is holomorphic, that is
$\deriv{u}{y}=J(u)\deriv{u}{x}$, to express all the deirvatives of $u$,
as derivatives with respect to $x$ only. We get~:
$$
\phi\circ u(z)=a_{2,0}x^{2}+a_{1,1}xy+a_{0,2}y^{2}+o(x^{2}+y^{2})
$$
where 
\begin{align*}
a_{2,0}&=\underset{a}{\underbrace{
         D^{2}\phi(X,X)
         }}
        +D\phi(\tsderiv[2]{u}{x^{2}})\\
a_{1,1}&=\underset{b}{\underbrace{
         D^{2}\phi(JX,X)
        +D\phi((X\cdot J)X)
        }}
        +D\phi(J\tsderiv[2]{u}{x^{2}})\\
a_{0,2}&=\underset{c}{\underbrace{
         D^{2}\phi(JX,JX)
        +D\phi(((JX)\cdot J)X-(X\cdot J)JX)
        }}
        -D\phi(\tsderiv[2]{u}{x^{2}})\\
\end{align*}
Recall from (\ref{eq:L(X)=h(X,X)+h(JX,JX)+new term}) that
$L_{\phi}(X)=D^{2}\phi(X,X)%
                +D^{2}\phi(JX,JX)%
                +D\phi(((JX)\cdot J)X-(X\cdot J)JX)$~:
\begin{equation}
  \label{eq:L(X)=a20+a02}
L_{\phi}(X)=a_{2,0}+a_{0,2}=\Delta(\phi\circ u)(0).
\end{equation}

Since $\deriv[2]{u}{x^{2}}(0)$ can be chosen arbitrarily, we can in
particular choose $d\phi(\deriv[2]{u}{x^{2}})(0)$ ({\it i.e.} its
$N$-component) and $d\phi(J\deriv[2]{u}{x^{2}})(0)$ ({\it i.e.} its $JN$
component) independently.

The choices%
\begin{align}
d\phi(\deriv[2]{u}{x^{2}})(0)&=\frac{c-a}{2}&
d\phi(J\deriv[2]{u}{x^{2}})(0)&=-b
\end{align}
lead to the formula
\begin{equation}
  \label{eq:phi o u=(a+c)/2(x2+y2)}
  \phi\circ u(z) = \frac{L_{\phi}(X)}{2}\ (x^{2}+y^{2}) + o(|z|^{2})
\end{equation}

If $L_{\phi}(X)>0$, $u$ is clearly tangent to $M$ from the outside.

If $L_{\phi}(X)=0$, then the contact order of $u$ with $M$ is at least
$3$.

If $u$ has a contact of order at least $3$ with $M$, then
$a_{2,0}=a_{1,1}=a_{0,2}=0$ and the relation (\ref{eq:L(X)=a20+a02})
proves that taking $X=\deriv{u}{x}(0)\in T^{J}_{0}M$, we have
$L_{\phi}(X)=0$.

This ends the proof of theorem \ref{thm:L>0 et courbe d'appui}.

\bigskip
The following definition of the higher order Levi forms derive from very
similar ideas in the study of the higher order Taylor expansion of
$\phi\circ u$~: the pairwise sums of the coefficients are related to the
derivatives of $\Delta(\phi\circ u)$, and thus to derivatives of the Levi
form.

Let $u$ be a germ of holomorphic disk at $0$. Using the relation
$\deriv{u}{y}=J(u)\deriv{u}{x}$ and eventually derivatives of it, it is
possible to express $\deriv[k]{\Delta(\phi\circ u)}{x^{p},y^{q}}$ by
means of derivatives of $u$ in the $x$ direction only, and this process is
purely algebraic~:
\begin{definition}\label{def:Lpq}
 For all $(p,q)$, there exists a polynomial function
 $L^{p,q}:(\R^{2n})^{p+q+1}\rightarrow\R$, called the $(p,q)$-th Levi
 form which depends on (the derivatives of) $\phi$ and $J$ only, such
 that for all germ of holomorphic disk $u$ at $0$~:
 $$
  \deriv[p+q]{}{x^{p},y^{q}}\ \Delta(\phi\circ u)
 =L^{p,q}(\deriv{u}{x},\dots,\deriv[p+q+1]{u}{x^{p+q+1}})
 $$
\end{definition}

At first sight, it seems that $L^{p,q}$ should depend on the derivatives
of $u$ up to order $p+q+2$ and not only $p+q+1$ as mentioned in this
definition. In fact, the derivatives of $u$ of order $p+q+2$ involoved in
the computation of $\deriv[p+q]{}{x^{p},y^{q}}\ \Delta(\phi\circ u)$
apear only in the following term~:
$d\phi(\deriv[p+q+2]{u}{x^{p+2},y^{q}}+\deriv[p+q+2]{u}{x^{p},y^{q+2}})$
which is null in the integrable case, and involves only derivatives of
$u$ of order at most $p+q+1$ in the non integrable one.

Especially in the almost complex case, the explicit computation of
$L^{p,q}$ is of course a bit tedious, but we stress that it is not more
complicated than the computation of the derivatives of some compound
function. More precisely, in the integrable case, $L^{p,q}$ can be
obtained by the following ``recipe''~: formally compute the $(p,q)$-th
derivative of $\Delta(\phi\circ u)$, and replace the
$(\alpha,\beta)$-derivatives of $u$ by $J^{\beta}X_{\alpha+\beta}$.

Another way to think about $L^{p,q}(X,X^{2},\dots,X^{k})$ for a given
vector field $X$, is to formally differentiate the usual Levi form, $p$
times in the $X$ direction, then $q$ times in the $JX$ direction, and
then ``force'' the commutation of $X$ and $JX$ by replacing all
$(JX)^{\beta}X^{\alpha}$ by $J^{\beta}X^{\alpha+\beta}$. In the non
integrable case, this commutation relation involves of course some
derivatives of $J$.

For instance, in the integrable case~:
\begin{align*}
  L^{0,0}(X_{1})
  &=D^{2}\phi(X_{1},X_{1})+D^{2}\phi(\i X_{1},\i X_{1})\\
  \begin{split}
    L^{1,0}(X_{1},X_{2})
    &=D^{3}\phi(X_{1},X_{1},X_{1})+D^{3}\phi(X_{1},\i X_{1},\i X_{1})+\\
    &\hspace{4cm}+2D^{2}\phi(X_{2},X_{1})+2D^{2}\phi(\i X_{2},\i X_{2})
  \end{split}\\
  \begin{split}
    L^{0,1}(X_{1},X_{2})
    &=D^{3}\phi(\i X_{1},X_{1},X_{1})+D^{3}\phi(\i X_{1},\i X_{1},\i X_{1})+\\
    &\hspace{4cm}+2D^{2}\phi(\i X_{2},X_{1})-2D^{2}\phi(X_{2},\i X_{2})
  \end{split}\\
  &\dots
\end{align*}

By construction, if $u$ is a \hideJ($J$-)holomorphic disk, the extended Levi
forms above are related to the Taylor expansion of $\phi\circ u$ at $0$ in
the following way~: write the homogeneous part of degree $k$ as 
$$
 [\phi\circ u(z)]_{k}=a_{k,0}x^{k}+a_{k-1,1}x^{k-1}y
+\dots+a_{0,k}y^{k}.
$$
Then
\begin{equation}\label{eq:Lij(du/dx)=ai+2,j + ai,j+2}
a_{i+2,j}+a_{i,j+2}=L^{i,j}(\tsderiv{u}{x},\dots,\tsderiv[k-1]{u}{x^{k-1}})
\end{equation}

We can now turn to the proof of theorem \ref{thm:Contact order via higher
  Levi forms}.

\subsection{Proof of theorem \ref{thm:Contact order via higher Levi
    forms}}

\paragraph{Proof of (\ref{item:Contact via Levi form-u
    exists})$\Rightarrow$(\ref{item:Contact via Levi form-Lpq(X)=0}) in
  theorem \ref{thm:Contact order via higher Levi forms}~:}

Start with a \hideJ($J$-)holomorphic disk $u$ tangent to $M$ at $0$ with
order $k+2$. Then the relation (\ref{eq:Lij(du/dx)=ai+2,j + ai,j+2})
proves that for all $(i,j), i+j+2\leq k+1$ we have
\begin{equation}
  \label{eq:Lpq(du/dx)=0}
L^{i,j}(\tsderiv{u}{x},\dots,\tsderiv[i+j+1]{u}{x^{i+j+1}})=0.  
\end{equation}
On the other hand, the corollary \ref{thm:u<->X} provides us with a
$J$-tangent vector field $X$ such that $\forall m\leq k+1$,
$\deriv[m]{u}{x^{m}}(0)=X^{m}(0)$, and therefore, $\forall(i,j), i+j\leq
k-1$~:
\begin{equation}
  \label{eq:Lpq(du/dx)=0}
L^{i,j}(X(0),\dots,X^{i+j+1}(0))=0.  
\end{equation}

\paragraph{Proof of %
 (\ref{item:Contact via Levi form-Lpq(X)=0})%
 $\Rightarrow$%
 (\ref{item:Contact via Levi form-u exists}) %
 in theorem \ref{thm:Contact order via higher Levi forms}~:}%
Start with a $J$-tangent vector field $X$ such that $\forall(i,j),
i+j\leq k-1,\, L^{i,j}(X(0),\dots,X^{i+j+1}(0))=0$. Choose a germ $u$ of
\hideJ($J$-)holomorphic disk at $0$ such that
\begin{equation}
  \label{eq:dudxm=Xm}
 \forall m\leq k+1,\ 
 \deriv[m]{u}{x^{m}}(0)=X^{m}(0).  
\end{equation}
Suppose by induction that $u$ is tangent to $M$ with order $k+1$. The
Taylor expansion of $\phi\circ u$ has the form
$$ 
 \phi\circ u(z)=a_{k+1,0}x^{k+1}+a_{k,1}x^{k}y
+\dots+a_{0,k+1}y^{k+1}+o(|z|^{k+1}).
$$
 and because of relation (\ref{eq:Lij(du/dx)=ai+2,j + ai,j+2}), we have 
\begin{equation}
  \label{eq:ai+2,j=-ai,j+2}
\forall (i,j), i+j+2=k+1,\ a_{i+2,j}=-a_{i,j+2}.
\end{equation}

Using (\ref{eq:dudxm=Xm}), one can easily compute $a_{k+1,0}$~:
\begin{align}
  a_{k+1,0}&=\deriv{}{x^{k}}\ d\phi\Big(\deriv{u}{x}\Big)\ (0)\notag\\
  &=X^{k}(d\phi(X))(0)=0.  \label{eq:ak+1,0=0}
\end{align}

To compute $a_{k,1}$, let us first remark that the computation of
$X^{k}\cdot(JX)(0)$ is a sum of terms involving a derivative of $J$ in
the $(X^{\alpha})_{\alpha\leq k+1}$ directions and one of the
$(X^{\alpha})_{\alpha\leq k+1}$. The computation of $\deriv{}{x^{k}}\
d\phi(J(u)\deriv{u}{x})\ (0)$ leads to the same sum where the
$X^{\alpha}$ are replaced by the $\deriv[\alpha]{u}{x^{\alpha}}$, and thus~:
\begin{align}
  a_{k,1}&=\deriv{}{x^{k}}\ d\phi\Big(J(u)\deriv{u}{x}\Big)\ (0)\notag\\
  &=X^{k}(d\phi(JX))(0)=0\label{eq:ak,1=0}
\end{align}

The relations (\ref{eq:ai+2,j=-ai,j+2}), (\ref{eq:ak+1,0=0}) and
(\ref{eq:ak,1=0}) prove that $\phi\circ u=O(z^{k+2})$, ending the proof
of the theorem.


\section{Base point dependency in lower dimensions}

In general, even in the pseudo-convex case, the type is not an upper
continuous function on $M$~: it is not even locally bounded (see
\cite{Da2} for striking examples). However, in complex dimension $2$, it
is known to be upper-continuous in the integrable case. In the first part
of this section, we give a proof of this fact with the point of vue
developped all along this paper~; in particular, the proof still apply to
the almost complex case.

In the second part of this section, we discuss a particular case in
complex dimension $3$~: if the Levi form has constant signature, the type
is an upper-continuous function of the base point.

\subsection{Upper continuity in dimension 2}

\begin{proposition}
In complex dimension $2$ (real dimension $4$), the regular type of an
hyper-surface $M$ is an upper-continuous function on $M$.
\end{proposition}

The basic idea of the proof is to translate $\typelisse(p)\geq k$ in
terms of jets~: the functions $L^{p,q}$ are defined on a space of jets,
and since they are invariant under rescaling, we can reduce to a compact
subset. The set $\{L^{p,q}=0,p+q\leq k\}$ is then compact and a sequence
of ``good'' jets above a converging sequence $p_{n}$ in $M$ should
therefore (sub-)converge to a ``good'' jet at $p_{\infty}$. The point is
that the set of realizable jets is not closed, since if the first
component of the jet vanishes, all the remaining components should vanish
as well.

This phenomenon occurs when the sequence of curves catches a singularity
in the limit. In terms of jets, if we rescale our jets so that the first
component has always length 1, this means that one component explodes in
the limit. Fortunately, in dimension $2$ (real dimension $4$), there is no
room for that, and the phenomenon can be controlled.

\begin{proof}

Let us first concentrate on one point where the type is $k+2$. Let $X$ be
a germ of $J$-tangent vector field at this point on which all the Levi
forms vanish up to order $k-1$.

Let $u$ be a germ of holomorphic disk at $0$ such that $\forall m\leq
k+1,\ \deriv[m]{u}{x^{m}}(0)=X^{m}(0)$. Then, $u$ has contact order $k+2$
with $M$ and thus, for any germ $\theta$ of holomorphic function with
$\theta(0)=0$, $\theta'(0)\neq 0$, the disk $v(z)=u\circ\theta(z)$ also
has contact order $k+2$. Suitably choosing the jet of $\theta$, and
computing the jet of $v$ at $0$, we obtain that in fact any
$Y=(\alpha+\beta J)X$ where $\alpha$ and $\beta$ are real valued
functions, {\it i.e.} any section of $T^{J}M$, satisfy the relations
$L^{p,q}(Y,\dots,Y^{p+q+1})=0$ for $p+q\leq k-1$.

Let now $X$ be a local non vanishing section $X$ of $T^{J}M$ near $0$,
and $p_{n}$ a sequence of points on $M$ converging to $0$. If each
$p_{n}$ is of type at least $k+2$, then the functions
$L^{p,q}(X,\dots,X^{p+q+1})$ vanish at each $p_{n}$, and thus at $0$.
therefore 
$$
\typelisse(\lim p_{n})\geq\limsup\typelisse(p_{n}).
$$

\end{proof}

\subsection{Case of dimension 3 with non degenerate Levi form}

The behaviour of the jet can also be controlled in some cases if the Levi
form is not degenerated~:
\begin{theorem}
  If $n=3$ and the Levi form on $M$ is (every where) non degenerate
  then $\typelisse(p)$ is an upper continuous function of $p$.
\end{theorem}

\begin{proof}
If the Levi form is definite positive (or negative) then
$\typelisse(p)=2$ every where.

If it is not, its signature is $(1,1)$, and the type is at least $3$
every where. Let $(X,Y)$ be a field of bases of $T^{J}M$ in which
$L=\big(\begin{smallmatrix}0&1\\1&0\end{smallmatrix}\big)$.
The coordinates are chosen such that $(X(0),Y(0),N(0))$ is the canonical
basis of $\C^{3}$. Suppose that the type is at least $k+2$ at $0$~: let
$u$ be a regular disk with contact $k+2$ and $Z$ a complex vector field on
$M$ realizing the jet of $\deriv{u}{x}$ up to order $k$.

In the integrable case, we can reparametrise $u$, so that
$u(z)=(z,h(z),g(z))$. In the non integrable case, this is not possible,
but the first coordinate of $u$ can still be brought to the form $z+r(z)$
where $r(z)=o(z)$ is a function whose derivatives are controlled by the
derivatives of $J$, and the derivatives of the other coordinates.

Let us decompose $Z^{k}(0)$ in the base $(X,JX,Y,JY,N,JN)$~:
$Z^{k}=\alpha X(0)+\beta Y(0) + \gamma N(0)$, where $\alpha,\beta,\gamma$ are
complex valued functions. 

The $X$-component of $Z^{k}$ is controlled by the $C^{k}$ norm of $J$ and
the other components of $Z^{k}$.
 
We already proved (see the proof of \ref{thm:Jets realisables}) that the
$N$ and $JN$ components of $Z^{k}$ can be expressed as a combination of
derivatives of $\phi$ (of length at most $k$) and all the $Z^{i}$ for
$i<k$~: we conclude that
  $$
  |\gamma|\leq \Vert\phi\Vert_{C^{k}} \sup_{i<k}\Vert Z^{i}\Vert^{k/i}
  $$

To control the $Y$ component of $Z^{k}$ at $0$, remark that in
$L^{k-1,0}(Z,\dots,Z^{k})$, the only term involving $Z^{k}$ is
$L(Z^{k},Z)$~:
$$
L^{k-1,0}(Z,\dots,Z^{k})=2L(Z^{k},Z)+P_{\phi}(Z,\dots,Z_{k-1})
$$
where $P_{\phi}(Z,\dots,Z_{k-1})$ is a homogeneous polynomial in
$Z,\dots,Z_{k-1}$ of degree $k+1$, whose coefficients come from
derivatives of $\phi$ and $J$ up to order $k$. 

We obtain that $L(\beta_{1}Y+\beta_{2}JY,X)=A_{\phi}(Z,\dots,Z_{k-1})$
where $A_{\phi}(Z,\dots,Z_{k-1})$ is an homogeneous polynomial in
$(Z,\dots,Z_{k-1})$ controlled by the $C^{k}$ norm of $\phi$ and $J$.

In the same way, we can prove that
$L(\beta_{1}Y+\beta_{2}JY,JX)=B_{\phi}(Z,\dots,Z_{k-1})$ where $B_{\phi}$
is again controlled by the $C^{k}$ norms of $\phi$ and $J$.

Finally
$$
\begin{cases}
  \beta_{1} L(Y,X) + \beta_{2}L(JY,X)&=A_{\phi}(Z,\dots,Z_{k-1})\\
  \beta_{1} L(Y,JX) + \beta_{2}L(JY,JX)&=B_{\phi}(Z,\dots,Z_{k-1})\\
\end{cases}
$$
and since $L(JY,X)=-L(JY,X)$, the determinant of this system is at least
$L(Y,X)^{2}=1$. As a consequence, $|\beta|$ is controlled by the $C^{k}$
norms of $\phi$ and $J$, and the norms of $(Z,\dots,Z_{k-1})$~: there
exists a contant $M$, depending only on the $C^{k}$ norms of $\phi$ and
$J$ such that
$$
\Vert Z^{k}\Vert\leq M \sup_{i<k}\Vert Z^{i}\Vert^{k/i}.
$$

By induction, we obtain a collection of constants $M_{i}$ such that at
$0$~:
$$
\forall i\leq k, 
\Vert Z^{i}\Vert\leq M_{i} \Vert Z\Vert^{i}.
$$

This inequality implies that if we consider a converging sequence $p_{n}$
of points of type at least $k+2$ in $M$, we can find a sequence of
realizable jets $(Z_{n},\dots,Z_{n}^{k})$, normalized so that $\Vert
Z_{n}(0)\Vert=1$, on which all the $L^{p,q}$ vanish for $p+q<k$, and such
that all the components stay bounded~: this sequence (sub-)converges to a
jet $\xi$, that is realisable, and on which all the $L^{p,q}$ vanish.
\end{proof}
\begin{remark}

We want to finish with an open question related to the base point
dependency of the type: in a complex case (see \cite{Da1}), we know that
the finite regular type is not an open condition ; but, if we consider
all curves (even singular), this condition become open. Is it true for a
almost complex structure and for pseudo-holomorphic curves ? We think
that the answer is positive.
\end{remark}


\begin{thebibliography}{99}
\bibitem{Blo}
T. Bloom.~: On the contact between complex manifolds and real hyper-surfaces in $\Bbb{C}^3$, Trans. Amer.
Math. Soc 263 (1981), 515-529.
\bibitem{Blo-Gra} 
T. Bloom-I. Graham.~: A geometric characterization of type on real submanifolds of $\Bbb{C}^n$, J. 
Diff. Geometry 12 (1977), 171-182.
\bibitem{Ca1}
D.Catlin.~: Subelliptic estimates for the $\bar\partial$-Neumann problem on pseudoconvex domains, Ann. of
Math. 126 (1987), 131-191.
\bibitem{Ca2}
D.Catlin.~: Estimates of invariant metrics on pseudoconvex domains of dimension two, Math. Z. 200 (1989),
429-466.
\bibitem{Da1}
J.-P. D'Angelo.~: Real hypersurface, orders of contact, and applications, Ann. of Math. 115 (1982),
615-637.
\bibitem{Da2}
J.-P. D'Angelo.~: Several complex variables and the geometry of real hypersurfaces, Studies in Adv. Math.
CRC press, Boca Raton, Fl, (1993).
\bibitem{Deb}
R. Debalme.~: Vari\'et\'es hyperboliques presque-complexes, Th\`ese de doctorat  de l'universit\'e de
Lille, (2001).
\bibitem{Kho}
J. Kohn.~: Boundary behavior of $\bar\partial$ on weakly pseudoconvexmanifolds of dimension two, J. Diff
Geometry 6 (1972), 523-542.
\bibitem{Mc}
J. Mcneal.~: Estimates on the Bergman kernels of convex domains, Adv. Math 109 (1994), 108-139.
\bibitem{Sik}
J.-C.~Sikorav.~: Some properties of holomorphic cueves in almost complex  manifolds, Prog. Math. 117,
165-189, Birkhauser, Basel, (1994).
\end{thebibliography}
\end{document}